\documentclass[12pt,twoside]{article}
\usepackage{amssymb, amsmath, amsthm,latexsym,graphics}
\usepackage[cp1251]{inputenc}
\usepackage[english]{babel}

\newtheorem{theorem}{Theorem}[section]
\newtheorem{lemma}[theorem]{Lemma}

\newtheorem{proposition}[theorem]{Proposition}

\newcommand{\seq}[2]{#1_1,#1_2,\ldots ,#1_#2}
\newcommand{\summa}[2]{\sum\limits_{#1}^{#2}}
\newcommand{\wt}[1]{\widetilde{#1}}
\newcommand{\mc}[1]{\mathcal{#1}}
\newcommand{\ov}[1]{\overline{#1}}
\newcommand{\src}[1]{\stackrel{\circ}{#1}}
\newcommand{\srb}[1]{\stackrel{\bullet}{#1}}
\newcommand{\RGU}{{\rm Rep}\,(G,\bigsqcup')}
\newcommand{\RGUC}{{\rm Rep}_\circ\,(G,\bigsqcup')}
\newcommand{\RGUB}{{\rm Rep}_\bullet\,(G,\bigsqcup')}

\begin{document}

\makeatletter
\renewcommand{\@biblabel}[1]{#1.}
\makeatother

\begin{center}
\LARGE \bf
Separating Functions, \\
Spectral Graph Theory \\
and Locally Scalar \\
Representations in Hilbert Spaces

\end{center}

\begin{center}
\Large\bfseries I.~K.~Redchuk~$^\dag$
\end{center}

\noindent
$^\dag$~Institute of Mathematics of National Academy of Sciences of Ukraine,\\
Tereshchenkovska str., 3, Kiev, Ukraine, ind. 01601\\
E-mail: red@imath.kiev.ua

\bigskip

\section{Separating Functions $\rho_r$}

The present paper is dedicated to studying the connections of separating
functions $\rho_r$, introduced in~\cite{RedRoi}, with locally scalar
representations of graphs on the one hand, and with spectral
graph theory on the other hand.

In the article~\cite{NazRoi} the function $P(S)$, attaching to
each partially ordered set $S$ a positive rational number,
was introduced. In terms of function $P$ the criterion of finite
representability and tameness of an arbitrary partially ordered set
is formulated in~\cite{NazRoi}: a patrially ordered set $S$ is finite
presented (tame) iff $P(S)<4$ ($P(S)=4$). If partially ordered
set $S$ is the union of $s$ disjoint chains (i.~e.
such that any two elements belonging to different chains are incomparable)
and $n_i$ is the number of elements in each chain, $i=\ov{1,s}$, then
$P(S)=\rho(\seq{n}{s})=\summa{i=1}{s}\rho(n_i)
=\summa{i=1}{s}1+\frac{n_i-1}{n_i+1},\; n_i\in\mathbb{N}.$ The list of all
solutions of the equations $\rho(\seq{n}{s})=4$ is:
\begin{equation}\label{eq11} (1,1,1,1),\;\; (2,2,2),\;\; (1,3,3),\;\;
(1,2,5).\;\; \end{equation} Since the function $\rho$ is increasing (the
partial order on integer vectors $(\seq{n}{s})$ is defined naturally), one
may easily obtain from the list~(\ref{eq11}) the list of all solutions of
the inequality $\rho(\seq{n}{s})<4$: \begin{equation}\label{eq11_}
(1,2,2),\;\; (1,2,3),\;\; (1,2,4),\;\; (1,1,k),\;\; (l,m),\;\; (n),\;\;
\end{equation} \noindent where $k,\, l,\, m,\, n$ are arbitrary natural
numbers.

It is remarkable that the list~(\ref{eq11}) corresponds to all
extended Dynkin graphs with one point of branching (the exact definition
see in s.~4 of the present paper) $\wt{D}_4$, $\wt{E}_6$, $\wt{E}_7$,
$\wt{E}_8$: the components of vectors from the list~(\ref{eq11})
correspond to the number of vertices on each branch. In the same way
the list~(\ref{eq11_}) describes Dynkin graphs $E_6$, $E_7$, $E_8$, $D_n$,
$A_n$.

For algebras defined by generators and linear or polylinear relations
(see~\cite{Red2}) the general formulas for dimensions can be obtained and
criterions of finite dimensionality and polynomiality of growth can be
also formulated in terms of the function $\rho$.

In~\cite{RedRoi} the natural generalization of the function
$\rho$ was suggested. Given $m\in\mathbb{N}$, denote $V_m$ the set
of unordered suits of $m$ nonnegative integers; then denote
$V = \bigcup\limits_{m\in\mathbb{N}} V_m$. Let for given
$r\in\mathbb{R}$, $r\ge 4$ a recurrent number sequence $\{u_i\}$ is
defined: $u_0=0$, $u_1=1$, $u_{i+2} = (r-2)u_{i+1}-u_i$. Put
\begin{equation}\label{eq12}
\begin{array}{ll}
\rho_r(0)=0,\\
\rho_r(n)=1+\frac{u_{n-1}}{u_n+1},\;\; n\in\mathbb{N},
\end{array}
\end{equation}
\noindent and for $\ov{v}\in V$, $\ov{v}=(\seq{n}{s})$
$$
\rho_r(\ov{v})=\summa{i=1}{s} \rho_r(n_i).
\footnote{The given definition slightly differs from one
in~\cite{RedRoi}: $\rho_r(n)$ here corresponds to $\rho_{r-4}$
in~\cite{RedRoi}.}
$$
It is easy to see that $\rho_4 \equiv \rho$.

The functions $\rho_r$ was exploited for describing
the standard characters of locally scalar representations of
certain types of graphs (see~\cite{RedRoi}, \cite{KruRoi}).

In~\cite{RedRoi} the following properties of the functions $\rho_r$ was
obtained:

\begin{proposition}\label{pro11}

For any $r > 4$
\begin{equation}\label{eq13}
\rho_r(n) = \frac{(\lambda+1)(\lambda^n-1)}{\lambda^{n+1}-1},
\end{equation}
\noindent where $\lambda = \frac{r-2+\sqrt{r^2-4r}}{2}$.

\end{proposition}

\begin{proposition}\label{pro12}

Given $r\in\mathbb{Z}$, the equation
\begin{equation}\label{eq14}
\rho_r(\seq{n}{s}) = r
\end{equation}
\noindent has following solutions:
$$
(\underbrace{1,1,\ldots,1}_{r}),\;\;
(\underbrace{2,2,\ldots,2}_{r-1}),\;\;
(\underbrace{1,3,3,\ldots,3}_{r-1}),\;\;
(\underbrace{1,2,5,5,\ldots,5}_{r-1}).
$$
For any $r\in\mathbb{Q}\setminus\mathbb{Z}$ the equation~(\ref{eq14}) has
no solutions.

\end{proposition}

Here we will give the simpler way of defining the functions $\rho_r$. Let
us show that \begin{equation}\label{eq15}
\rho_r(n+1)=\frac{r}{r-\rho_r(n)}
\end{equation}
\noindent for any $n\in\mathbb{N}\bigcup\{0\}$.

Indeed, by formula~(\ref{eq13}) we have $\frac{r}{r-\rho_r(n)} =
\frac{r}{r-\frac{(\lambda+1)(\lambda^n-1)}{\lambda^{n+1}-1}}$; since
$r=\frac{(\lambda+1)^2}{\lambda}$ and $\lambda\neq 1$, the last
equality reduces to form
$\frac{(\lambda+1)(\lambda^n-1)}{(\lambda+1)(\lambda^{n+1}-1) -
\lambda(\lambda^n-1)} =
\frac{(\lambda+1)(\lambda^{n+1}-1)}{\lambda^{n+2}-1} = \rho_r(n+1)$, which
required.

\medskip

Further we will define the functions $\rho_r$ by the formulas~(\ref{eq15})
for $r\ge 1$, $r\in\mathbb{R}$ (the starting condition is the same:
$\rho(0)=0$.) However for $1\le r < 4$ we may obtain
$\lambda^{n+1}=1$ and $\rho_r(n)$ will be not defined. The equality
$\lambda^{n+1}=1$ is equivalent to $\lambda = \cos \frac{2\pi
k}{n+1} + i\sin \frac{2\pi k}{n+1}$, $k=\ov{1,n}$ ($\lambda\neq 1$), which
means that $r = 4\cos^2 \frac{\pi k}{n+1}$, $k=\ov{1,n}$. When
$r=4\cos^2\frac{\pi k}{n+1}$, $k=\ov{1,n}$ we will formally
put $\rho_r(n) = \infty$, and all rational transformations with
such $\rho_r(n)$ we will do using the natural transition to the limit. In
particular, if $\rho_r(n) = \infty$ then $\rho_r(n+1) = 0$, which
implies that for $r=4\cos^2\frac{\pi k}{n+1}$, $k=\ov{1,n}$, the function
$\rho_r(m)$ is periodic with period $n+1$.

It is easy to determine that the formulas~(\ref{eq12}) and (\ref{eq13}) and
the proposition~\ref{pro12} hold for such definition.

Let us point to one more case, where the functions $\rho_r$
appear naturally. $*$-representations of $*$-algebras $\mc{P}_{r,\alpha} =
\mathbb{C}\langle \seq{p}{r}\, |\, p_k^2 = p_k^* = p_k,\;
\summa{k=1}{r}=\alpha e \rangle$ in $H$, where $e$ is the identity of
algebra, $H$ is separable Hilbert space, are studied in~\cite{KruRabSam}.
One of the results of this work is such that if
$\Sigma_r$ is the set of those $\alpha\in\mathbb{R}$ for which
$\mc{P}_{r,\alpha}$ has at least one representation and $r\ge 4$, then
$\Sigma_r = \{\Lambda_r^1,\, \Lambda_r^2,\, \left[
\frac{r-\sqrt{r^2-4r}}{2},\, \frac{r+\sqrt{r^2-4r}}{2},\, r-\Lambda_r^1,\,
r-\Lambda_r^2 \right]\}$, where $\Lambda_r^1$, $\Lambda_r^2$ are discrete
sets, which can be defined recurrently:
\begin{eqnarray}
\Lambda_r^1 =
\{0,1+\frac{1}{r-1},1+\frac{1}{(r-2)-\frac{1}{r-1}},\ldots,\nonumber\\
1+\frac{1}{(r-2)-\frac{1}{(r-2)-\frac{1}{\ddots_{-\frac{1}{r-1}}}} } \}
\nonumber
\end{eqnarray}
\begin{eqnarray}
\Lambda_r^2 =
\{1,1+\frac{1}{r-2},1+\frac{1}{(r-2)-\frac{1}{r-2}},\ldots,\nonumber\\
1+\frac{1}{(r-2)-\frac{1}{(r-2)-\frac{1}{\ddots_{-\frac{1}{r-2}}}}} \},
\nonumber
\end{eqnarray}
It is easy to see, that
$\Lambda_r^1=\{\rho_r(2k)\}$, $\Lambda_r^2=\{\rho_r(2k+1)\}$,
i.~e.
$\Sigma_r= \{\{\rho_r(k)\}, \left[\frac{r-\sqrt{r^2-4r}}{2},\,
\frac{r+\sqrt{r^2-4r}}{2},\, r-\Lambda_r^1,\, r-\Lambda_r^2 \right],
\{r-\rho_r(k)\}\}$, $k\in\mathbb{N}\bigcup\{0\}$\footnote{The introduction
of the functions $\rho_r$ was proposed by A.~V.~Roiter exactly for defining
of these sets.}.

In~\cite{KruRabSam} on the categories $\mathrm{Rep}\,\mc{P}_{r,\alpha}$
of $*$-representations of algebras $\mc{P}_{r,\alpha}$ it was determined
functors $\Phi^+$ and $\Phi^-$ (Coxeter functors), which gave an
opportunity to describe all irreducible $*$-representations of
algebras $\mc{P}_{r,\alpha}$ (up to unitary equvalence) in points
of discrete spectrum of the set $\Sigma_r$. There was also given the explicit
(but rather complicated) formula for calculating
$\Phi^{+k}(\alpha)=(\Phi^+)^k(\alpha)$. Let us show that this formula in
terms of the functions $\rho_r$ has more simple form. To do that,
we will need

\begin{lemma}\label{lem13}

$\rho_r(2k-1) = 1+\frac{u_{k-1}}{u_k}.$

\end{lemma}

{\bf Proof.} If $r=4$ then $\rho_r(2k-1) = 1+\frac{2k-2}{2k} =
1+\frac{k-1}{k} = 1+\frac{u_{k-1}}{u_k}$.

If $r>4$ then $u_k = \frac{\lambda^k - \lambda^{-k}}{\sqrt{r^2-4r}}=$
(see~\cite{RedRoi}) $=\frac{\lambda^{2k}-1}{\lambda^{k-1}(\lambda^2-1)}$.
Therefore, $1+\frac{u_{k-1}}{u_k} = 1+\frac{\lambda^{2k-2}-1}{\lambda^{k-2}}
\cdot \frac{\lambda^{k-1}}{\lambda^{2k}-1} = \rho_r(2k-1)$. The lemma is
proved.

\medskip

We will show now that
$$
\Phi^{+k}(\alpha) =
\frac{r-\rho_r(2k-1)\alpha}{r-\rho_r(2k-1)-\alpha}.
$$

In~\cite{KruRabSam} it was obtained that $\Phi^{+k}(\alpha) =
1+\frac{a_{k-1}}{a_k}$, where $a_k$ is defined recurrently: $a_1 = 1$, $a_2
= r-1-\alpha$, $a_{k+2} = (r-2)a_{k+1} - a_k$. Simple induction gives $a_k
= u_{k+1} + (1-\alpha)u_k$. Thus, $\Phi^{+k}(\alpha) =
1+\frac{a_{k-1}}{a_k} = 1+\frac{u_k+u_{k-1}-\alpha u_{k-1}}{u_{k+1}+u_k -
\alpha u_k} \overbrace{=}^{(\ref{eq12})}
\frac{ru_k-\alpha(u_k+u_{k-1})}{ru_k-(u_k+u_{k-1})} =
\frac{r-\rho_r(2k-1)\alpha}{r-\rho_r(2k-1)-\alpha}$, which required.

\section{Spectra and indexes of graphs}

Further (if the contrary is not indicated specially) all
considering graphs suppose to be finite, connected and not containig
loops and multiple edges.

Let $G_v$ to be the set of vertices and $G_e$ be the set of edges of graph
$G$. Determine the numeration on the vetrices of graph $G$: $G_v =
\{\seq{g}{n}\}$, $n\in\mathbb{N}$. Denote $M(g_k)=\{g_i\, |\, g_i
\mbox{ connected with } g_k\}$, $1\le k,i\le n$.

Matrix $A_G = ||a_{ij}||_{i,j = \ov{1,n}}$, where
$a_{ij} =
\begin{cases}
1, & \text{if $g_i\in M(g_j)$,}\\
0, & \text{if $g_i\not\in M(g_j)$}
\end{cases}$
\noindent is called the {\it adjacency matrix} of graph $G$.

Since matrix $A_G$ is symmetrical, all its eigenvalues over
$\mathbb{C}$ are real. The linearly ordered set $\sigma(G) =
\{\lambda_{min} = \lambda_1\le \lambda_2 \le \ldots \le \lambda_n =
\lambda_{max}\}$ of eigenvalues of the matrix $A_G$ is called the {\it
spectrum} of graph $G$, and the number $\mathrm{ind}(G) =
\lambda_{max}$ is the {\it index} of the graph $G$.

Denote $V_G$~--- linear space over $\mathbb{R}$,
consisting of suits $x=(x_i)$, identifying each vector $x\in V_G$ with
function $x: G_v\to \mathbb{R}$, $x_i = x(g_i)$. Elements $x\in V_G$
are called {\it $G$-vectors} (see~\cite{KruRoi}). The adjacency matrix $A_G$
can be considered as the matrix of certain linear oparator in the
natural basis in the space $V_G$.

The spectral graph theory has been studied rather deeply
(see~\cite{CveDoobZac}, \cite{CveRowSim}). In particular, the following
remarkable statement, known as Smith's theorem, holds:

\begin{theorem}\label{thS}\cite{Smi}

Let $\lambda = \mathrm{ind}(G)$ for connected graph $G$. Then $\lambda <
2$ iff $G$ is Dynkin graph ($A_n$, $D_n$, $E_6$,
$E_7$, $E_8$); $\lambda = 2$ iff $G$ is
extended Dynkin graph ($\wt{A}_n$, $\wt{D}_n$, $\wt{E}_6$,
$\wt{E}_7$, $\wt{E}_8$).

\end{theorem}

Note, that if $G$ is disconnected, then its index is equal to the maximal
of the indexes of its connected components. Also the following statements
hold (\cite{CveDoobZac}, \cite{CveRowSim}):

\begin{proposition}\label{pro21}

For an arbitrary graph $G$ the inequalities $1\le \mathrm{ind}(G)
\le |G_v|-1$ hold.

\end{proposition}

\begin{proposition}\label{pro22}

The index $\lambda = \mathrm{ind}(G)$ of an arbitrary graph
$G$ is a simple eigenvalue, iff graph $G$ is connected,
and in this case the eigenspace of $V_G$,
belonging to $\lambda$, is spanned by a vector, whose coordinates
are all positive.

\end{proposition}

This vector is called the {\it principle eigenvector} of $G$. If
$\lambda = \mathrm{ind}(G)$, $A_G = ||a_{ij}||$, then the condition
``$x=(\seq{x}{n})$ is a principal eigenvector of $G$'' can be
written in a form
\begin{equation}\label{eq21}
\lambda x_i = \summa{j=1}{n}a_{ij}x_j,\;\; i,j = \ov{1,n}.
\end{equation}
graph $G$ is called {\it bipartite}, if $G_v = \srb{G}_v\bigsqcup
\src{G}_v$, $\srb{G}_v \bigcap \src{G} = \varnothing$ and $g \in
\srb{G}_v$ implies $M(g)\subseteq \src{G}_v$ and conversely: $h \in
\src{G}_v$ implies $M(g)\subseteq \srb{G}_v$. The set $\src{G}_v$
is called the set of {\it even verices} and the set
$\srb{G}_v$ is the set of {\it odd vertices} of $G$
(see~\cite{KruRoi}.) If the graph $G$ is bipartite and at first
we numerate add, and then~--- even vertices, then its adjacency
matrix has the form
$$
A_G = \begin{bmatrix}
O_1 & B \\ B^* & O_2
\end{bmatrix}
$$
\noindent where $O_1$, $O_2$ are quadratic zero matrices
of the orders $|\srb{G}_v|$ and $|\src{G}_v|$ respectively, $B^*$
is a transposed matrix of $B$.

Evidently, any tree is bipartite graph, and a cycle with $n$ vertices
is bipartite iff $n$ is even. (Therefore,
all Dynkin graphs and extended Dynkin graphs are bipartite,
except $\wt{A}_{n-1}$ with odd number of vertices $n$.)

\section{Standard characters of star-shaped graphs}

Let $G$ be a bipartite graph, $G_v = \srb{G}_v\bigsqcup \src{G}_v$.
Determine numeration of vertices of $G$ such that $\srb{G}_v =
\{\seq{g}{p}\}$, $\src{G}_v = \{g_{p+1},g_{p+2},\ldots,g_n\}$. Let
(according to this numeration) $y = (\seq{y}{n})$ is a principal
eigenvector of $G$, i.~e. vector $y$ satisfies the equalities~(\ref{eq21}):
\begin{equation}\label{eq31}
\lambda y_i =
\summa{j=1}{n}a_{ij}y_j,\;\; i,j = \ov{1,n},
\end{equation}
\noindent where $\lambda = \mathrm{ind}(G)$.
Vector $y_\bullet = (\seq{y}{p},0,0,\ldots,0)$ is
the {\it odd standard vector}, and vector $y_\circ =
(0,0,\ldots,0,y_{p+1},y_{p+2},\ldots,y_n)$ is the {\it even standard
vector} of $G$. Since the principal eigenvector $y$ is determined
up to nonzero real multiplier, $y_\bullet$ and
$y_\circ$ are also determined up to common for both vectors
$y_\bullet$, $y_\circ$ nonzero multiplier.

The {\it Coxeter transformation} in the space $V_G$ is a linear
transformation $c= \sigma_{g_n}\cdot\cdots\cdot\sigma_{g_1}$, where
$\left(\sigma_{g_i}(x)\right)_j = x_j$ for $i\neq j$ and
$\left(\sigma_{g_i}(x)\right)_i = -x_i + \sum\limits_{j\,|\,g_j\in
M(g_i)} x_j$. Denote $\srb{c} =
\sigma_{g_p}\cdot\cdots\cdot\sigma_{g_1}$, $\src{c} =
\sigma_{g_n}\cdot\cdots\cdot\sigma_{g_{p+1}}$, i.~e. $c=\src{c}\srb{c}$.
Obviously, $c^{-1} = \srb{c}\src{c}$ and $(\srb{c})^2 = (\src{c})^2 =
\mathrm{id}$. Put $c_t = \underbrace{\cdots \srb{c}\src{c}\srb{c}}_{t}$
for $t > 0$, $c_t = \underbrace{\cdots \src{c}\srb{c}\src{c}}_{t}$ for $t <
0$ and $c_0 = \mathrm{id}$. For $G$-vectors $v$ and $w$ we denote
$v\simeq w$, if $v$ and $w$ are linearly dependent.

In~\cite{RedRoi} the explicit formulas in terms of the functions
$\rho_r$, indicating how odd and even
standard vectors transform under the action of $c_t$, was
obtained for certain graphs of the special type. Now we
will prove an analogous statement for an arbitrary bipartite graph $G$.

\begin{proposition}\label{pro31}

Let $G$ be a bipartite graph, $r=(\mathrm{ind}(G))^2\ge 4$, $y_\bullet$,
$y_\circ$ is its odd and even standard vectors,
$t\in\mathbb{N}$. Then the following hold:
\begin{equation}\label{eq32}
\begin{array}{llll}
c_{2t-1}(y_\circ) \simeq
\left(y_\circ +
\frac{\sqrt{r}}{\rho_r(2t-1)}y_\bullet\right);\\
c_{2t}(y_\circ) \simeq
\left(y_\circ +
\frac{\rho_r(2t)}{\sqrt{r}}y_\bullet\right);\\
c_{-(2t+1)}(y_\bullet) \simeq
\left(y_\bullet +
\frac{\sqrt{r}}{\rho_r(2t+1)}y_\circ\right);\\
c_{-2t}(y_\bullet) \simeq
\left(y_\bullet +
\frac{\rho_r(2t)}{\sqrt{r}}y_\circ\right).\\
\end{array}
\end{equation}

\end{proposition}

{\bf Proof.} We will prove the statement by induction for the case
$c_m(y_\bullet)$, where $m<0$ (in the case $m>0$ the proof is exactly
the same up to ``evenness''.)

For $t=-1$ we have $c_{-1}(y_\bullet) = \src{c}(y_\bullet) =
\src{c}(y_1,\ldots,y_p,0,\ldots,0) = (y_1,\ldots,y_p,x_{p+1},\ldots,x_n)$,
where $x_k = \sum\limits_{g_i\in M(g_k)} y_i = \sqrt{r}y_k$, $k=\ov{1,p}$
(formulas (\ref{eq31})). Thus, $c_{-2}(y_\bullet) = (r-1)y_\bullet +
\sqrt{r}y_\circ \simeq y_\bullet + \frac{\sqrt{r}}{r-1}y_\circ = y_\bullet
+ \frac{\rho_r(2)}{\sqrt{r}}y_\circ$, since $\rho_r(2) = \frac{r}{r-1}$.

Assume, that the formulas~(\ref{eq32}) hold for all $m\le 2t$.
Then $c_{-(2m+1)}(y_\bullet) = \src{c}c_{-2t}(y_\bullet) \simeq
\src{c}(y_\bullet + \frac{\rho_r(2t)}{\sqrt{r}}y_\circ) =
\src{c}(y_1,\ldots,y_p,\frac{\rho_r(2t)}{\sqrt{r}}y_{p+1},\ldots,\frac{\rho_
r(2t)}{\sqrt{r}}y_n) = \src{c}(y_1,\ldots,y_p,x_{p+1},\ldots,x_n)$, where
$x_k = \sum\limits_{g_i \in M(g_k)}y_i - \frac{\rho_r(2t)}{\sqrt{r}}y_k =
y_k(\sqrt{r}-\frac{\rho_r(2t)}{\sqrt{r}}) =
\frac{r-\rho_r(2t)}{\sqrt{r}}y_k = \frac{\sqrt{r}}{\rho_r(2t+1)}y_k$,
therefore, $c_{-(2m+1)}(y_\bullet) = y_\bullet +
\frac{\sqrt{r}}{\rho_r(2t+1)}y_\circ$.

Then, $c_{-(2t+2)}(y_\bullet) = \srb{c}c_{-(2t+1)}(y_\bullet) =
\srb{c}(y_\bullet + \frac{\sqrt{r}}{\rho_r(2t+1)}y_\circ) =
\srb{c}(y_1,\ldots,y_p,\frac{\sqrt{r}}{\rho_r(2t+1)}y_{p+1},\ldots,
\frac{\sqrt{r}}{\rho_r(2t+1)y_n}) = (x_1,\ldots,x_p,
\frac{\sqrt{r}{\rho_r(2t+1)}y_{p+1},\ldots,
\frac{\sqrt{r}}{\rho_r(2t+1)}y_n})$, where $x_k =
\frac{\sqrt{r}}{\rho_r(2t+1)} \sum\limits_{g_i\in M(g_k)} y_i - y_k =
y_k(\frac{r}{\rho_r(2t+1)}-1)$. Therefore, $c_{-(2t+2)}(y_\bullet) =
(\frac{r}{\rho_r(2t+1)}-1)y_\bullet + \frac{r}{\rho_r(2t+1)}y_\circ \simeq
y_\bullet + \frac{\sqrt{r}}{\rho_r(2t+1)} \cdot
\frac{\rho_r(2t+1)}{r-\rho_r(2t+1)}y_\circ = y_\bullet + \frac{1}{\sqrt{r}}
\cdot \frac{r}{r-\rho_r(2t+1)}y_\circ = y_\bullet +
\frac{\rho_r(2t)}{\sqrt{r}}y_\circ$.

The proposition is proved.\footnote{These formulas were also obtained
independently by V.~L.~Ostrovskyi ({\tt arxiv: math.RA/0509240})
in the case of star-shaped graphs.}

\medskip

Consider now the connection between standard vectors of a bipartite
graph $G$ and locally scalar representations of this graph. Locally
scalar representations of graphs were introduced and studied
in~\cite{KruRoi}. We remind some notions from this work.

Let $\mc{H}$ be the category of Hilbert spaces, which objects are
separable Hilbert spaces, and morphisms are bounded
operators. {\it Representation $\pi$ of a graph} $G$ in $\mc{H}$
attaches to each vertex $a\in G_v$ an object
$\pi(a)=H_a\in\mathrm{Ob}\,\mc{H}$ and to each edge $\gamma \in
G_e$ connecting vertices $a$ and $b$ a pair of interadjoint linear
operators $\pi(\gamma) = \{\Gamma_{ab},\Gamma_{ba}\}$, where $\Gamma_{ab}:
H_b\to H_a$. Denote $A_g = \sum\limits_{b\in M(g)}
\Gamma_{gb}\Gamma_{bg}$, $b,g \in G_v$. A representation $\pi$ is called
{\it locally scalar}, if all operators $A_g$ are scalar, $A_g =
\alpha_gI_{H_g}$, where $I_{H_g}$ is identity operator in space
$H_g$. Since operators $A_g$ are bounded, $\alpha_g \ge 0$. $G$-vector
$f$, such that $f(g) = \alpha_g$ for all $g\in G_v$, is called the {\it
character} of the representation, and $G$-vector $d$, such that $d(g)=\dim
\pi(g)$ is the dimension of the representation $\pi$. Note, that for
given representation $\pi$ its dimension is determined uniquely, but
the character, in general, is not (it is determined uniquely on the {\it
support} $G^\pi =\{a\in G_v\,|\, \pi(a)\neq 0\}$ of the representation.)
A representation $\pi$ is {\it faithful}, if $G^\pi = G_v$.

In~\cite{KruRoi} for an arbitrary graph $G$ such categories are considered:
a category $\rm{Rep}\,(G,\mathcal{H})$ of representations of $G$ in the
category of Hilbert spaces, a category $\rm{Rep}\,(G)$ of locally
scalar representations, its (full) subcategory $\rm{Rep}\,(G,d,f)$
of representations with fixed dimension $d$ and character $f$.
In~\cite{RedRoi} for bipartite graph $G$ are considered: a category
$\RGU$~--- the union of categories $\rm{Rep}\,(G,d,f)$, such that $f(g)>0$
for all $g\in M(G^\pi)$; a category
$\RGUC \subset \RGU$~--- full subcategory of representations with
character $f$, for which $(d,f)$, such that  $f(g)>0$ for $g\in
(M(G^\pi)\bigcup G^\pi)\bigcap \srb{G}_v$, and an analogous category
$\RGUB$; also there was defined functor $\src{\Phi}$, which is an
equivalence of the category $\RGUC$, such that $\src{\Phi}(f)_i = f_i$ for
$g_i \in \src{G}_v$, and if $g_i \in \srb{G}_v$, then $\src{\Phi}(f)_i =
f_i$ for $d_i = 0$, $g_i \in M(G^\pi)$, and in another cases
$\src{\Phi}(f)_i = \sigma_i(f)$; and the analogous functor $\srb{\Phi}$ is
defined. Also for any $k\in\mathbb{N}$ functors $\Phi_k=\underbrace{\cdots
\srb{\Phi}\src{\Phi}\srb{\Phi}}_{k}$ and $\Phi_{-k}=\underbrace{\cdots
\src{\Phi}\srb{\Phi}\src{\Phi}}_{k}$ are constructed, such that
$d(\Phi_t(\pi)) = c_t(d(\pi))$, $t\in\mathbb{Z}$.

A nonnegative $G$-vector $x$ is called {\it regular}, if $c_t(x)$ is
nonnegative for any $t\in\mathbb{Z}$, and {\it singular} in the opposite
case. Locally scalar representation $\pi$ of a graph $G$ is {\it singular}
({\it regular}), if $\pi$ is indecomposable, and $d(\pi)$ is
a singular (regular) vector. An object $(\pi,f)\in\mathrm{Ob}\,\RGU$ is {\it
singular}, if $\pi$ is singular.

A group $W$, generated by reflections $\sigma_{g_i}$, is called the {\it
Weyl group}. A vector $x\in V_G$ is called the ({\it real)
root}, if $x= w\ov{a}$ for certain $a\in G_v$ and $w\in W$, where
$\ov{a}$ is a {\it simple root}, i.~e. $\ov{a}(a)=1$ and $\ov{a}(a) = 0$
for $g\neq a$.

The {\it simplest object} in the category $\RGU$ is a pair
$(\Pi_g,\ov{f})$, such that $\dim \Pi_g = \ov{g}$; $\ov{f}(g) = 0$ and
$\ov{f}(a) > 0$ for $a\in M(g)$.

\begin{theorem}\label{th32}\cite{RedRoi}

Each singular object of the category $\RGU$ can be obtained as
$\Phi_m(\Pi_g,\bar{f})$, where $(\Pi_g,\bar{f})$ is a
simplest object ($m \geq 0$ for $g \in \src{G}_v$ and $m \leq 0$ for $g \in
\srb{G}_v$). At that each faithful singular representation $G$
corresponds (up to equivalence) to one singular object of $\RGU$.

\end{theorem}

A character $f$ of a representation is {\it standard}, is $f =
c_t(y_\circ)$ or $f=c_t(y_\bullet)$, $t\in\mathbb{Z}$. A representation
$(\pi, f)$ with a standard character $f$ is a {\it standard representation},
and the corresponding object of the category $\RGU$ is a {\it standard
object.}

The following statement is well-known (see., for instance, \cite{Kac})

\begin{lemma}\label{lemRoot}

If $x$ is a real root, ther either $x > 0$, or $(-x) > 0$.

\end{lemma}

Let us prove the following

\begin{proposition}\label{proSingRoot}

If $d$ is a positive real singular root of a graph $G$, then
there exists such $t \in \mathbb{Z}$ that $d = c_t(\ov{g})$, $g\in G_v$.

\end{proposition}

{\bf Proof.} Let $d$ be not a simple root (in the opposite case
we have $t=0$.) Then, let $m$ be a minimal in absolute value integer
number with property $c_m(d) < 0$. Let, for definiteness, $m >
0$. Then $c_{m-1}(d) > 0$ (lemma~\ref{lemRoot}). Therefore,
$c_{m-1}(d)$ is  a simle root (if $c_{m-1}(d)$ has at least two
positive coordinates, then it is clear, that coordinates, corresponding
to the neighbour vertices, are zeros, and, applying the reflection in one
of these positive coordinates, we obtain a contradiction with
lemma~\ref{lemRoot}.) Thus, $t=1-m$ and the proposition is proved.

\begin{theorem}\label{th33}

Let $G$ be a bipartite graph with $\mathrm{ind}(G)\ge 2$, $d$ is a
singular real root in $G$. Then there exists a unique standard
representation $\pi$ with dimension $d$.

\end{theorem}

{\bf Proof.} Existence. Since $d$ is a singular real root, there exists $t\in
\mathbb{Z}$, such that $d = c_t(\ov{g})$ (proposition~\ref{proSingRoot}).
Therefore, allowing theorem~\ref{th32} and proposition~\ref{pro31}, we have
for $t\ge 0$ $\Phi_t(\Pi_g,y_\circ)$ and for $t\le 0$
$\Phi_t(\Pi_g,y_\bullet)$ is equal to $(\pi,y)$, where $y$ is a standard
character, $\dim \pi = d$.

Uniqueness. Formula (\ref{eq32}) imply that if $(\pi,f)$ is a
standard object of the category $\RGU$, then $\src{\Phi}(\pi,f)$ is a
standard object only if $\pi\neq \Pi_g$, where $g\in\src{G}_v$, and
$\srb{\Phi}(\pi,f)$ is a standard object, only if $\pi\neq \Pi_g$,
where $g\in\srb{G}_v$.

\section{Indexes of star-shaped graphs}

A {\it path} of length $l\in\mathbb{N}$ on a graph $G$ is an ordered
sequence of vertices $(g_{i_1},\ldots,g_{i_{l+1}})$,
such that $g_{i_k}\in M(g_{i_{k+1}})$, $k=\ov{1,l}$. Vertices $g_{i_1}$ and
$g_{i_{l+1}}$ are called the beginning and the end of path
respectively. A vertex $g\in G_v$ of $G$ is called the {\it point
of branching}, if $M(g)\ge 3$. {\it Star-shaped
graph} is a tree which has no more than one point of branching.
For a star-shaped graph $G$ the set $G_v$ can be presented as
$$
G_v = B_0\bigsqcup B_1\bigsqcup B_2 \bigsqcup \ldots \bigsqcup B_s,
$$
\noindent where $B_i\bigcap B_j = \varnothing$, $i,j=\ov{0,s}$,
$B_0=\{g_0\}$, and $g_0$ is a point of branching in $G$ (if one
exists) and $g,h\in B_i$ for certain $i$ iff the path of
minimal length with the beginning in $g$ and the end in $h$ does
not contain $g_0$. The sets $B_i$ are {\it branches} of $G$.

Further we will point at the direct connection of the separating functions
$\rho_r$ with indexes of star-shaped graphs.

First we need to prove the following lemma.

\begin{lemma}\label{lem41}

Given real $r\ge 1$. Let $\{v_n\}$ be a number sequence, defined
recurrently: $v_0=0$, $v_1=1$, $v_{n+2} =
\sqrt{r}v_{n+1}-v_n$. Then $\rho_r(n) = \sqrt{r}\frac{v_n}{v_{n+1}}$.

\end{lemma}

{\bf Proof.} Induction by $n$. $\rho_r(0) =
\frac{\sqrt{r}v_0}{v_1} = 0$. Assume $\rho_r(n) =
\frac{\sqrt{r}v_n}{v_{n+1}}$. Then by formula (\ref{eq15}) $\rho_r(n+1) =
\frac{r}{r-\frac{\sqrt{r}v_n}{v_{n+1}}} =
\frac{\sqrt{r}v_{n+1}}{\sqrt{r}v_{n+1}-v_n} =
\frac{\sqrt{r}v_{n+1}}{v_{n+2}}$, which required.

\medskip

\begin{theorem}

Let $G_v = B_0\bigsqcup B_1\bigsqcup \ldots \bigsqcup B_s$ be the
separation of a star-shaped graph by branches, $|B_i|=n_i$, $i=\ov{1,s}$.
Then $\rho_r(\seq{n}{s})=r$, where $r=(\mathrm{ind}(G))^2$.

\end{theorem}

{\bf Proof.} Let $B_0 = \{g_0\}$, $B_k =
\{g_1^k,\ldots,g_{n_k}^k\}$, $k=\ov{1,s}$, and vertices in the branches $B_k$
are numerated in such a way: $M(g_1^k)=1$, $g_i^k\in M(g_{i+1}^k)$,
$i=\ov{1,n_k-1}$, $g_{n_k}^k\in M(g_0)$. Let $y$ be a principal eigenvector
of $G$, $y_i^k = y(g_i^k)$, $y_0=y(g_0)$, $k=\ov{1,s}$, $i=\ov{1,n_k}$.
Then the equations (\ref{eq21}) have the form \begin{equation}\label{eq41}
\begin{array}{lllll} \lambda y_1^k = y_2^k,\\ \lambda y_2^k =
y_1^k+y_3^k,\\ \lambda y_3^k = y_2^k+y_4^k,\\ \cdots\; \cdots\; \cdots \\
\lambda y_{n_k}^k = y_{n_{k-1}}^k+y_0,\ \end{array} \end{equation}
\begin{equation}\label{eq42} \lambda y_0 = y_{n_1}^k+\cdots +y_{n_k}^k,
\end{equation} where $\lambda = \mathrm{ind}(G)$.

The equations (\ref{eq41}) for any branch $B_k$ imply $y_0 =
v_{n_k+1}y_1^k$ and $y_{n_k}^k = v_{n_k}y_1^k$, where $\{v_i\}$ is
a sequence, defined recurrently: $v_0 = 0$, $v_1=1$,
$v_{i+2}=\lambda v_{i+1}-v_i$. Then $y_{n_k}^k =
\frac{v_{n_k}}{v_{n_k+1}}y_0$ for all $k=\ov{1,s}$. Consider
these relations and the equation~(\ref{eq42}) we obtain $\lambda y_0 =
\summa{k=1}{s} \frac{v_{n_k}}{v_{n_k+1}} y_0$ and, taking into account
that $y_0 \neq 0$ (proposition~\ref{pro22}) and $\lambda = \sqrt{r}$,
we have $\sqrt{r} = \summa{k=1}{s}\frac{v_{n_k}}{v_{n_k+1}}$. Then, by
lemma~\ref{lem41} $r=\summa{k=1}{s}\rho_r(n_k)=\rho_r(\seq{n}{s})$, which
required.

\end{document}